\documentclass[12pt]{amsart}
\usepackage{graphicx}
\usepackage{mathptmx}
\usepackage{latexsym,amsmath,amssymb,amsfonts,amsthm}

\newcommand{\sm}[1]{\mbox{\small $#1$}}

\newcommand{\La}[1]{\mbox{\Large $#1$}}

\newcommand{\ha}{\sm{\frac{1}{2}}}

\newcommand{\R}{\mathbb{R}}          

\newfont{\lie}{eufm10 at 12pt}
\newfont{\field}{msbm10 at 11pt}
\newtheorem{theorem}{Theorem}[section]
\newtheorem{lemma}{Lemma}[section]
\newtheorem{corollary}{Corollary}[section]
\newtheorem{proposition}{Proposition}[section]
\newtheorem{remark}{Remark}
\title{A spectral Bernstein theorem}
\author{Pedro  Freitas \and Isabel Salavessa}
\date{\today}
\thanks{Partially supported by 
FCT through project PTDC/MAT/101007/2008.}
\address{Department of Mathematics, Faculdade de Motricidade Humana (TU Lisbon)
{\rm and} Group of Mathematical Physics of the University of Lisbon\\ Complexo
Interdisciplinar, Av. Prof. Gama Pinto 2\\ P-1649-003 Lisboa,
Portugal}\email{freitas@cii.fc.ul.pt}
\address{Centro de F\'{\i}sica das Interac\c{c}\~{o}es
Fundamentais, Instituto Superior T\'{e}cnico, Technical University
of Lisbon, Edif\'{\i}cio Ci\^{e}ncia, Piso 3, Av.\ Rovisco Pais,
P-1049-001 Lisboa, Portugal}\email{isabel.salavessa@ist.utl.pt}
\begin{document}
\maketitle
{\small {\bf Abstract:}  We study the spectrum of the Laplace operator
of a complete
minimal properly immersed hypersurface $M$ in $\R^{n+1}$.
(1) Under a volume growth condition on extrinsic balls and a condition
on the unit normal at infinity, we prove that $M$ has only essential spectrum
consisting of the half line $[0, +\infty)$. This is the case when 
$\lim_{\tilde{r}\to +\infty}\tilde{r}\kappa_i=0$, where  $\tilde{r}$ 
is the extrinsic distance to a point of $M$ and $\kappa_i$ are the 
principal curvatures.
(2) If the $\kappa_i$  satisfy the decay conditions
$|\kappa_i|\leq 1/\tilde{r}$,
and strict inequality is achieved at some point $y\in M$, 
then there are no eigenvalues. 
We apply these results to minimal graphic and multigraphic hypersurfaces.}
\section{Introduction}
The graphic minimal equation in $\R^{n+1}$, for a
function $f:\R^n\to \R$, is given by
$$\sum_i\left(\frac{D_if}{\sqrt{1+|Df|^2}}\right)=0.$$
It is well known that  entire solutions of this equation are
 linear  if $n\leq 7$ (see \cite{B,Fl,dG,Al,Si}), and there are 
counterexamples for $n\geq 8$ given by Bombieri, De Giorgi and  Giusti
\cite{BdGG}. A natural question to ask is whether these submanifolds may be
distinguished by their spectral properties or not.

The Laplace operator $-\Delta$ on a complete noncompact 
Riemannian manifold $M$ acting on $C^{\infty}_0(M)$ is essentially 
self-adjoint and can be uniquely extended as an unbounded self-adjoint
operator to a subspace $\mathcal{D}$ of $L^2(M)$ of functions $u$ for which
$\Delta u\in L^2(M)$ in the sense of distributions.
In what follows, by spectrum of $M$ we mean the spectrum of $-\Delta $.
This is  a closed subset of $ [0, +\infty)$, which can be decomposed as
$\sigma(M)=\sigma_p(M)\cup \sigma_{ess}(M)$, where
$\sigma_p(M)$ is the pure point spectrum, composed by
isolated eigenvalues of finite multiplicity,
and $\sigma_{ess}(M)$ is the essential
spectrum, a closed subset of $[0, +\infty)$, which
is the set of values of $\lambda$ for which there exists
an $L^2$-orthonormal  sequence $u_m$ 
with $ (\Delta +\lambda I)u_m \to 0$, in $L^2(M)$.
The essential spectrum also includes eigenvalues of infinite multiplicity,
 and those of finite multiplicity that are not isolated, if they exist.
The bottom of the spectrum admits a variational formulation and may be
determined via the minimization of the Rayleigh
quotient, namely,
$$\lambda(M)= \inf_{u\in C^{\infty}_0(M)}\frac{
\int_M\|\nabla u\|^2dV}{\int_Mu^2dV}.$$
In the case where $M$ has infinite volume, the lowest point of the essential
spectrum $$\lambda_{ess}=\inf \sigma_{ess}(M)\geq \lambda(M),$$
can be estimated by Brooks's inequality \cite{Br}
$$\lambda_{ess}\leq \frac{1}{4}\mu^2_M,$$
where $\mu_M$ is the exponential volume growth of $M$
$$ \mu_M=\lim\sup_{r\to +\infty} \frac{1}{r}\log(V_M(B_r(x))),$$
and where $V_M(B_r(x))$ is the volume of the geodesic ball in $M$
of radius $r$ and center $x$ ( $\mu_M$ does not depend on $x$).
If the Ricci tensor of $M$ is bounded from below by a constant $K\in \R$
then, by the Bishop volume comparison theorem, $\mu_M$ is finite
(and in particular $\sigma_{ess}(M)$ is nonempty), and is zero if $K\geq 0$.
In general, minimal submanifolds of $\R^{n+k}$ do not have to satisfy this 
property, for 
it corresponds to bounded second fundamental form. 

Recall that the Euclidean space $\R^n$ has only essential spectrum consisting of
the whole half line $[0,+\infty)$ and there are no (embedded) eigenvalues. 
We may ask which minimal submanifolds of $\R^{n+k}$ have a trivial spectrum,
that is, which minimal submanifolds have the same spectrum as 
$\R^n$. In the case where a minimal submanifold is properly immersed
in a ball  of $\R^{n+k}$, then it is known that there exists only
pure point spectrum, as proved by Bessa, Jorge and Montenegro~\cite{BJM}.
The situation for unbounded submanifolds will, however, be different in general,
and one may ask this question for particular cases such as minimal graphs in
$\R^{n+1}$. A first step 
towards answering this
is to determine whether or not there exist minimal graphic hypersurfaces with the same
trivial spectrum as $\R^{n}$ as described above. In this paper we present
some results towards an answer to this question -- see Theorems~\ref{thm1}
and~\ref{thm2} below. As an example, we recover a result that may already
be obtained from the
work of Donnelly~\cite{Do2}, that the catenoid surface in $\mathbb{R}^{3}$ has 
spectrum $[0, \infty)$, and further prove -- Corollary~\ref{catenoid} -- 
that it has no embedded
eigenvalues. Thus minimal multigraphic hypersurfaces cannot be distinguished by their spectra. It
remains open whether or not this is also the case for minimal graphic hypersurfaces.

In what follows, $F:M\to \R^{n+1}$
is a complete, oriented, properly immersed
minimal hypersurface. On $M$ we give the induced metric $g_M$ and
corresponding Levi-Civita connection $\nabla$, 
 and  denote by $A\in C^{\infty}(\odot^2TM^*\otimes NM)$ 
the second fundamental form of $F$, $A(X,Y)=D_XY -\nabla_XY$, where $D$ 
stands for the flat connection on
$\R^{n+1}$, and $NM$ is the normal bundle of $M$.
We denote by $\tilde{r}$ the distance
function in $\R^{n+1}$ to a fixed point $F(x)$, and by
$\tilde{B}_r(x)\cap M$,  the pull back  by $F$ of  
the ball $\tilde{B}_r(F(x))$ on $\R^{n+1}$, with center
$F(x)$  and radius $r$, and  call it the extrinsic ball at $x\in M$. 
It contains the intrinsic ball $B_r(x)$.
\begin{theorem} \label{thm1}Assume $M$ is a complete oriented properly immersed minimal
 hypersurface of
$\R^{n+1}$, $x\in M$, and that there exists a positive constant
$C_n>0$ such that, for all $r>0$ sufficiently large
\begin{equation}\label{cond1}
 V_M(\tilde{B}_r(x)\cap M)\leq C_nr^n, 
\end{equation}
where $V_M$ is the volume with respect to the induced metric $g_M$ of $M$.
Then $\mu_M=0$. Furthermore, if 
\begin{equation}\label{cond2}
|d\tilde{r}(\nu)|\to \xi,~~~~~~~\mbox{when}~\tilde{r}\to +\infty, 
\end{equation}
where $0\leq \xi< 1$ is a constant,
and $\nu$ is the unit normal to $M$, then $\sigma(M)=\sigma_{ess}(M)
=[0,+\infty)$.
\end{theorem}
Condition~(\ref{cond1}) implies the number $\kappa(M)$
 of ends of $M$
must be less than or equal to $C_n/\omega_n$ where $\omega_n$ is the
volume of the unit ball in $\R^n$ (see \cite{Chen}). 
The first part of next proposition is obtained from Lemma 2.4
and Theorem 2.2 of the work of Q.\ Chen 
\cite{Chen} (see Remark 3 in Section 3):
\begin{proposition} If $\lim_{\tilde{r}(F(y))\to +\infty}
\tilde{r}(F(y))\|A(y)\|=0$,
then $(\ref{cond1})$ is satisfied with $C_n=\kappa(M)\omega_n$, and 
$(\ref{cond2})$ is also satisfied with $\xi=0$, for $n\geq 3$ or,
for $n=2$, provided $M$ has embedded ends.
In particular, $\sigma(M)=\sigma_{ess}(M)=[0,+\infty)$.
\end{proposition}
\noindent
Note that the decay condition on $A$ in this
proposition is not satisfied by the 
examples of non-linear minimal graphs of Bombieri, de Giorgi and Giusti.

The next theorem is well known for stable varifolds
(see for instance \cite{Simon}, Theorem 17.7~\footnote{We are indebted
to Brian White for calling our attention to this result.})
and shows a reverse inequality to~(\ref{cond1}):
\begin{theorem} (Volume monotonicity formula)
 Let $M$ be a properly immersed minimal hypersurface
of $\R^{n+1}$. 
Then for each $x\in M$
and 
for any $r>\epsilon>0$ $$\frac{V_M(\tilde{B}_r(x)\cap M)}{r^n}
\geq \frac{  V_M(\tilde{B}_{\epsilon}(x)\cap M)}{\epsilon^n}
=:F_n(x,\epsilon).$$
Furthermore, $\lim_{\epsilon\to 0}F_n(x,\epsilon) =k\omega_n$,
where $k$ is the number of self-intersections of $F$ at $x$, and
$\omega_n$ is the
volume of  the unit ball of $\R^n$.
\end{theorem}
Given a function  $f:\R^n\to \R$,  we denote the graphic hypersurface
by   $M= \Gamma_f=\{(p,f(p)):p\in \R^n\}$. 
If $\Gamma_f$ is minimal then inequality  (1) is 
satisfied by a classical result due to Miranda:
\begin{theorem}[\cite{Mi}] If $\Gamma_f$ is an entire minimal graphic 
hypersurface, then, for each $p\in \R^n$
and $r>0$ 
$$\int_{\{q\in \R^n: (q,f(q))\in \tilde{B}_r(p,f(p))\cap \Gamma_f\}}
(1+|Df|^2)^{1/2}dV\leq \frac{(n+1)^2}{2}\omega_{n+1}r^n,$$
where  $dV$ is the Euclidean volume element of $\R^{n}$. Furthermore
$$\int_{\{q\in \R^n: (q,f(q))\in \tilde{B}_r(p,f(p))\cap \Gamma_f\}}
\|A\|^2dV_M\leq k(n)r^{n-2},$$
where $k(n)$ is a constant depending on $n$.
\end{theorem}
\noindent
The left-hand-side of the first inequality 
is just the volume of the extrinsic ball 
$V_M(\tilde{B}_r(x)\cap M)$.
Thus, we obtain as a corollary of the main theorem:
\begin{corollary} If $M=\Gamma_f$ is a minimal graphic hypersurface,
 then $\mu_M=0$. 
Furthermore, if  there is $x=(p,f(p))$
satisfying~ $(\ref{cond2})$, then
$\sigma(M)=\sigma_{ess}(M)=[0,+\infty)$
\end{corollary}
Condition~(\ref{cond2}) holds in several different situations.
For instance, we will see 
in Lemma~\ref{lemma1} that if $M=\Gamma_f$ for a function
 $f:\R^n\to \R$,
and if there exists
a unit vector $p_0\in \R^n$ such that $\gamma(t)=f(tp_0+p)$ has
bounded derivative, then
$
\lim\inf_{\tilde{r}\to +\infty}$ 
$ |d\tilde{r}(\nu)|=0.
$
If $|Df|$ is bounded  and $\Gamma_f$ 
is minimal, then
 Moser in (\cite{Mo}) proved that for each $k$, $u=\frac{\partial f}
{\partial x_k}$ is solution of a uniformly elliptic 
second order differential 
equation  in the selfadjoint form. Deriving a suitable  Harnack  theorem
that allows for the estimation of the growth of the oscilation of $u$,
he concludes that $\lim_{|q|\to +\infty} Df_q$ exists,
and further that $f$ is a linear affine function.
 In fact from
the existence of the previous limit, we  can show that 
(\ref{cond2}) holds with
$\xi=0$ (see Lemma 3.1). 
Linear maps $f(q)=l(q)+b$, satisfy $d\tilde{r}(\nu)=0$.
 More generally, (\ref{cond2}) holds
with $\xi=0$ 
 if  at infinity, for $q$ within open sets, 
$f(q)$ is   of the form 
$C|q|^{2\alpha}$,  $C\log|q|$, or $Ce^{B|q|^{2\alpha}}$,
where $C,B, \alpha$
are any reals. On the other hand, since $\Gamma_f$ is a minimal
graph, then $f$ is a harmonic function for the graph metric, which implies
that $f$ satisfies maximum and minimum principles. In particular, $f$
cannot be globally of the form $\phi(|q|)$ in some ball $B^0_s(0)$ of $\R^n$, 
for some fixed function $\phi$,  unless it is constant.
For a given minimal graphic hypersurface, there could exist
more than one limit point of $|d\tilde{r}(\nu)|$, when $\tilde{r}\to
+\infty$, where $\nu$ is the unit normal to the hypersurface.
It is not clear to us at this point if this may have implications on the 
spectral behaviour.

To prove Theorem~\ref{thm1}, for each $\lambda>0$ we build
a sequence $u_{m}$ spanning an infinite dimensional subspace of
$L^2(M)$ and such that $(\Delta+\lambda)u_{m}\to 0$.
This is achieved  by using test functions
supported in annuli of extrinsic balls and using the volume growth
estimates.
This is a similar construction given by J.\ Li \cite{Li} where
intrinsic balls were used.

Next we give a  decay condition on the second fundamental form $A$ of $M$
 that implies the non-existence of eigenvalues. 
\begin{theorem}\label{thm2} If $F:M\to \mathbb{R}^{n+1}$ is a complete
properly immersed minimal hypersurface such that the second fundamental
form satisfies $\|A(X,X)\|\leq {1}/{\tilde{r}}$, 
for any unit tangent vector $X$, and strict inequality is achieved
at some point of $M$, 
then $M$ has no eigenvalues.
\end{theorem}
\begin{remark}{\rm 
If $M$ has finite total scalar curvature, that is $\int_M\|A\|^ndV_M<+\infty$
(this condition is sufficient to ensure a  complete minimal
immersed submanifold in a Euclidean space is properly immersed,
see  \cite{Anderson}), then Q.\ Chen in \cite{Chen} proved that 
$\lim_{r\to +\infty} V_M(\tilde{B}_r(x)\cap M)/\omega_nr^n$ is just
$\kappa(M)$. We also note that, under the assumption of finite
total scalar curvature,
Anderson \cite{Anderson} concluded that
$\|A\|\leq c/\tilde{r}^n$, for some constant $c$, and that
for $r$ sufficiently large one has
$\sup_{\partial \tilde{B}_r(x)}\|A\|\leq \mu(r)/r$
where $\mu(r)\to 0$ when $r\to +\infty$. 
These  inequalities with respect to $\|A\|$ also show that
the assumption in Theorem~\ref{thm2} above is quite natural.}
\end{remark}

\noindent
The condition on $A$ is equivalent to a similar condition on the
 principal curvatures $\kappa_i$ of $M$. The proof of theorem 1.4 consists on
a similar  construction in~\cite{EF}, 
using now the extrinsic distance function instead of the intrinsic
one.
\begin{corollary} If $\Gamma_f$ is a minimal graphic hypersurface
defined by a function $f:\mathbb{R}^n\to \mathbb{R}$ such that
$f(0)=0~$ and
$$|D^2f(x)(X,X)|^2 \leq (1+ |Df(x)|^2)/
(|x|^2+ f^2(x)), $$
for any $|X|\leq 1$, with strict inequality at some $x$, 
 then $\Gamma_f$ has no eigenvalues. Furthermore, if $(2)$ is satisfied,
then $\Gamma_f$ has trivial spectrum.
\end{corollary}
For a multigraph we have the example of the catenoid. It is not difficult to
see that conditions (1)  and (2) of Theorem 1.1
are satisfied, and  that the principal curvatures $\kappa_i$ 
satisfy  the conditions in Proposition 1.1 and in Theorem 1.4.
To prove all this we only have to
recall that $t\leq \sinh t\cosh t$ and $\sinh t\leq t\cosh t$,  
for any $t\geq 0$,  with equality only at $t=0$, as we
can verify using an infinite Taylor expansion of the hyperbolic sine
and cosine. 

\begin{corollary}\label{catenoid} The  catenoid  surface in $\R^3$ has
 trivial spectrum.
\end{corollary}
Minimal hypersurfaces of Euclidean spaces have nonpositive
Ricci tensor, and  by  the generalization of
the Hilbert-Efimov theorem given by Smith and Xavier \cite{SX},  
 $\inf \|A\|=\sup \mathrm{Ricci}^M=0$. No
further information on the curvature is given.
These results are insufficient to allow us to apply known results relating
the spectrum (of a minimal hypersurface) to curvature 
(see e.g \cite{ChenZ,Do1,Do2,DoPL}), for the case $n> 2$.
For $n=2$, Donnelly in Theorem 6.3 of \cite{Do2} proved that 
complete noncompact simply connected 
surfaces with nonpositive  curvature that converges to zero at infinity, 
have  essential spectrum $[0,+\infty)$. For surfaces with finite fundamental 
group he obtained only $\lambda_{ess}=0$
with no further conclusion on the nonexistence of eigenvalues.
This last result may be applied to the catenoid, but 
Corollary~\ref{catenoid} above characterizes the whole
spectrum. We note that,
for an arbitrary noncompact
complete Riemannian manifold,  it is sufficient that
$Ricci^M$ converges to zero at infinity in order to  have $\lambda_{ess}= 0$.
This can be derived from an argument used in \cite{Do2},
applying   Cheng's
eigenvalue comparison  inequality \cite{Cheng}, as was used in
\cite{CEK} (see also Proposition 3.2).

\section{Proof of Theorem 1.2}
We will give a proof of Theorem 1.2  with no need of using
currents, closely following the proof of a recent result
due to Alencar, Walcy and Zhou:
\begin{theorem}[\cite{AWZ}]  If $M$ is a minimal immersed hypersurface
of $\R^{n+1}$ with induced metric $g_M$,
then for each $x\in M$ and $r>0$, away from cut locus distance of $x$,
$$V_M(B_r(x))\geq \omega_n r^n.$$
\end{theorem}
\noindent
By using extrinsic balls, we avoid  the assumption
on the cut locus. We denote by $\tilde{g}$ the Euclidean metric
of $\R^{n+1}$, and by $g_M$ the induced metric on  $M$. For $z\in \R^{n+1}$,
 $\tilde{r}(z)=|z-F(x)|=\sqrt{\tilde{g}(z-F(x),z-F(x))}$.
Restricting $\tilde{r}$ to $M$, $\tilde{r}(y)=\tilde{r}(F(y))$, we set
$h=\sm{\frac{1}{2}}\tilde{r}^2:M\to \mathbb{R}$.
A standard computation shows that, for $X\in T_yM$, 
$d\tilde{r}_y(X)=\tilde{g}(X, F(y)-F(x))/\tilde{r}(y)$ and 
\begin{equation}
Hess\, h_y(X,X)=\tilde{g}(X,X)+ \tilde{g}(A(X,X),F(y)-F(x)).
\end{equation}
Therefore, if $M$ is a minimal hypersurface with mean curvature $H$, $nH=\mathrm{tr}\, A=0$,
and so,
$$ \Delta h=n.$$
Integration of the previous equation
 on a normal domain $D$ of $M$ with boundary $\partial D$
  and Stokes's Theorem yields
\begin{equation}\label{stokes}
\int_{\partial D}\tilde{r}\, g_M(\nabla\tilde{r},\nu_{\partial D})\, dS
= nV_M(D),
\end{equation}
where $\nu_{\partial D}$ is the unit normal to  ${\partial D}$ and 
 $dS$ is the volume element, and $\nabla\tilde{r}$ the gradient in $M$. 
Taking $D=\tilde{B}_s(x)\cap M$,
where $s$ is a regular value of $\tilde{r}_{|M}$, define
$$
V_M(s):=V_M(\tilde{B}_s(x)\cap M).
$$
Using~(\ref{stokes})  with 
$\nu_{\partial \tilde{B}_s}=\nabla\tilde{r}/\|\nabla\tilde{r}\|$,
and the co-area formula,
\begin{equation}
\int_{\partial \tilde{B}_s(x)\cap M}\tilde{r}\, 
\|\nabla \tilde{r}\| dS
= nV_M(s)=n\int_0^s\frac{1}{t}\rho(t) dt,
\end{equation}
where ~$
\rho(t)=\int_{\partial\tilde{B}_t(x)\cap M}
\frac{\tilde{r}}{\|\nabla\tilde{r}\|}dS.~$ ~
Since $\|\nabla\tilde{r}\|^{-1}\geq 1
\geq \|\nabla\tilde{r}\|$, by (5), 
$$\rho(s)\geq \int_{\partial\tilde{B}_s(x)\cap M}
\tilde{r} \|\nabla\tilde{r}\|dS = n \int_0^s\frac{1}{t}\rho(t)dt.
$$
Therefore $s\, V_M'(s) \geq n V_M(s)$ which implies
$\frac{d}{ds}\ln V_M(s)\geq \frac{d}{ds} \ln s^n $. Integration along
$[\epsilon,s] $, where $0<\epsilon<s$, leads to
$
{V_M(s)}/{s^n}\geq {V_M(\epsilon)}/{\epsilon^n}, 
$
and the inequality of the theorem is proved.

Next we prove $\lim_{\epsilon\to 0}V_M(\epsilon)/\epsilon^n=k\omega_n$.
We take ${\gamma}(t)$ a curve in $M$
starting at $x$. Thus,  $\tilde{\gamma}(t):=F(\gamma(t))=
F(x)+t\tilde{\gamma}'(0)+o(t)$,
with $\tilde{\gamma}'(0)$ non zero. Since
$\tilde{\nabla}\tilde{r}(\tilde{\gamma}(t))=\tilde{\gamma}(t)-F(x)/
|\tilde{\gamma}(t)-F(x)|$, then
$$
\tilde{g}(\tilde{\nabla}\tilde{r}(\tilde{\gamma}(t)),\nu_{\gamma(t)})
=\frac{t\tilde{g}(\tilde{\gamma}'(0), \nu_{\gamma(t)})+o(t)}{
|t\tilde{\gamma}'(0)+o(t)|}
= \frac{\tilde{g}(\tilde{\gamma}'(0), \nu_{\gamma(t)})+\frac{o(t)}{t}}
{|\tilde{\gamma}'(0)+\frac{o(t)}{t}|}
$$
and this converges to $0$ when $t\to 0$.
For $y\in M$, $1=|\tilde{\nabla}\tilde{r}|^2=\|\nabla \tilde{r}\|^2
+ |d\tilde{r}(\nu)|^2$, and so ~
$\lim_{y\to x, \, y \in M}\|\nabla \tilde{r}\|
=1$. Therefore,
\begin{equation}
\lim_{\epsilon\to 0}\frac{1}{V_S(\partial\tilde{B}_{\epsilon}(x)\cap M)}
\int_{\partial \tilde{B}_{\epsilon}(x)\cap M}
\frac{1}{\|\nabla \tilde{r}\|}dS=1.
\end{equation}
Making $\epsilon\to 0$, using l'H\^{o}pital's Rule, (5) and (6), we have
\begin{eqnarray*}
\lim_{\epsilon\to 0} \frac{V_M(\epsilon)}{\epsilon^n}&=&
\lim_{\epsilon\to 0}\frac{\rho(\epsilon)}{n\epsilon^n}
= \lim_{\epsilon\to 0}\frac{1}{n\epsilon^{n-1}}
\int_{\partial\tilde{B}_{\epsilon}(x)\cap M}
\frac{1}{\|\nabla\tilde{r}\|}dS\\
&=& \lim_{\epsilon\to 0}\frac{V_S(\partial\tilde{B}_{\epsilon}(x)\cap M)}
{n\epsilon^{n-1}}=
\lim_{\epsilon\to 0}\frac{V_S(\partial\tilde{B}_{\epsilon}(x)\cap M)}
{V_{n-1}(S^{n-1}_{\epsilon})}\omega_n=k\omega_n,
\end{eqnarray*}
where $k$ is the  number of self-intersections of $F(M)$ at $F(x)$, 
and $S^{n-1}_{\epsilon}$ is the sphere of $\R^n$ of radius $\epsilon$. 
In the last equality we used that,  when $\epsilon\to 0$, 
$F(M)\cap \tilde{B}_{\epsilon}(F(x))$ can be identified with
$k$ copies of
$T_xM=\R^n$. \qed
\section{Proof of the Theorem 1.1 and Corollary 1.1}
Since $|F(x)-F(y)|\leq d(x,y)$ for any $x,y\in M$,  where $d$ is the 
intrinsic distance on $M$, then
$B_r(x)\subset \tilde{B}_r(x)\cap M$. In particular
\begin{equation}\label{l7}
\mu_M\leq \lim\sup_{r\to +\infty}\frac{1}{r}\ln\left(\int_{\tilde{B}_r(x)
\cap M}dV_M\right).
\end{equation}

For the case  $M=\Gamma_f$,
we denote by $g$ and 
$\tilde{g}=g+dt^2$ the Euclidean metrics of $\R^n$ and 
$\R^{n+1}=\R^n\times\R$, respectively.
We often identify $\Gamma_f$, a subset of $\R^{n+1}$, with
$\R^n$ endowed with the graph metric
$g_M= g+ f^*dt^2$,
that is, we may see $\Gamma_f$ as an immersion
$\Gamma_f:\R^n\to \R^{n+1}$,
$\Gamma_f(p)=(p,f(p))$ and
give to $\R^n$ the pull back  metric by $\Gamma_f$.
The graph metric is always  complete.
In this case, 
$B_r^0(p)$ denotes a ball for the Euclidean metric $g$ while
$B_r(p)$ is a ball with respect to the graph metric $g_M$.
 The volume element of
$\R^n$ with respect to the graph metric is given by
\begin{equation}\label{l8}
 dV_M=\sqrt{1+|Df|^2}dV, 
\end{equation}
where $Df$  denotes the $g$-gradient of $f$, and the unit normal to the 
graph is 
\begin{equation}
\nu= \frac{(-Df,1)}{\sqrt{1+|Df|^2}}. 
\end{equation}
We have, $B_{r'}^0(p)\subset B_r(p)$, where $r=r'\sqrt{1+F(p,r')^2}~$
with $F(p,r')=$ $\sup_{q\in B_{r'}(p)}|Df|$, and 
$B_r(p)\subset $ $\tilde{B}_r(p,f(p))\cap \Gamma_f\subset B_r^0(p).$ 

From~(\ref{l7}),~(\ref{l8}) and  Miranda's Theorem 1.3, 
we conclude:
\begin{proposition} 
If $\Gamma_f$ is minimal, then it has zero exponential volume 
growth, $\mu_M=0$.
In particular $0\in \sigma_{ess}(M)$, and it is not an eigenvalue.
\end{proposition}
\noindent
\em Proof. \em 
 By a result due to Yau \cite{Y}, there are no
$L^2$ harmonic functions on noncompact complete manifolds
of infinite volume (see Theorem 1.2). \qed
\begin{remark}{\rm From Miranda's inequality and the volume monotonicity
formula we see that $\frac{\omega_{n+1}}{\omega_n}
\geq \frac{2}{(n+1)^2}$. There is a sharper lower bound,
$\frac{\omega_{n+1}}{\omega_n}
\geq \frac{\sqrt{2\pi}}{\sqrt{n+2}}$, \cite{Bor}.}
\end{remark}
\begin{lemma}\label{lemma1} $(i)$~ If $M$ is a minimal hypersurface of $\R^{n+1}$, 
then at any $y\in M$
$$
\frac{n-1}{\tilde{r}}\leq  \Delta \tilde{r}\leq \frac{n}{\tilde{r}}, \quad
\quad 1-\|\nabla \tilde{r}\|^2 =
\tilde{g}(\tilde{\nabla}\tilde{r},\nu)^2.
$$
$(ii)$~In the case $M=\Gamma_f$, $x=(p,f(p))$, at $y=(q,f(q))$, 
\begin{equation}
\tilde{g}(\tilde{\nabla}\tilde{r},\nu)=
\frac{1}{\tilde{r}}\frac{(-g(q-p,Df_q)+(f(q)-f(p))}{\sqrt{1+|Df|^2}}.
\end{equation}
Moreover, if there exists a direction $p_0\in \R^n$, $|p_0|=1$,
such that $\gamma(t)=f(tp_0+p)$ is a curve such that 
$\exists\lim_{t\to +\infty}\gamma'(t)$ and is bounded  then
$$
\lim_{t\to +\infty} \tilde{g}(\tilde{\nabla}\tilde{r}_{\Gamma_f(\gamma(t))},
\nu_{\Gamma_f(\gamma(t))})^2~=0.$$
In particular, if $\exists \lim_{\tilde{r}\to +\infty}Df$ then
$(\ref{cond2})$ is satisfied with $\xi=0$.
\end{lemma}
\noindent
\em Proof. \em  Let $e_i$ be a  $g_M$ o.n. basis of $T_yM$. We denote by
 $(u)^{\top}$ and $(u)^{\bot} $ the orthogonal projections
of a vector $u\in \R^{n+1}$  onto $T_yM$
and  $NM_y$, respectively.
For any $y\in M$, from $|\tilde{\nabla}\tilde{r}|^2=1$ and that
$(\tilde{\nabla}\tilde{r})^{\top}=\nabla \tilde{r}$,
$(\tilde{\nabla}\tilde{r})^{\bot}=\tilde{g}(\tilde{\nabla}\tilde{r},\nu)\nu$, 
we obtain the second equality. 
Since the mean curvature $H$ of $M$ vanishes,
a standard computation gives
 at $y$
\begin{eqnarray*}
\Delta\tilde{r}&=& \sum_i \tilde{D}d \tilde{r}\, (dF(e_i),dF(e_i))
+\tilde{g}(nH, \tilde{\nabla}\tilde{r})=  
\sum_i \tilde{D}d \tilde{r}\, (dF(e_i),dF(e_i)) \\
&=& \frac{n}{\tilde{r}}-\sum_i
\frac{\tilde{g}(F(y)-F(x),e_i)^2}{\tilde{r}^3}
= \frac{n}{\tilde{r}}-\frac{|(F(y)-F(x))^{\top}|^2}{\tilde{r}^3}, 
\end{eqnarray*}
Now 
${|(F(y)-F(x))^{\top}|^2}\leq 
{|F(y)-F(x)|^2}\leq
\tilde{r}^2$, 
 and we obtain the bounds of $\Delta\tilde{r}$.
For $M=\Gamma_f$, $x=(p,f(p))$,
$y=(q,f(q))$, (9) is expressed as  (10). Now we 
assume $\exists\lim_{t\to +\infty}\gamma'(t)$ and is bounded.
Let $q(t)=tp_0+p$.
Along
$(q(t), \gamma(t)=f(q(t)))$,  by (10),
for all $t$ sufficiently large,
\begin{eqnarray*}
|\tilde{g}(\tilde{\nabla}\tilde{r},\nu)| 
&=& \frac{1}{\tilde{r}}\frac{\La{|}(-t\gamma'(t)+\gamma(t)-f(p))\La{|}}
{\sqrt{1+|Df_{q(t)}|^2}}.
\end{eqnarray*}
Therefore, multiplying by $1/t$, we have
\begin{eqnarray*}
\lefteqn{\La{|}(-\gamma'(t)+\frac{\gamma(t)}{t}-\frac{f(p)}{t})\La{|}
~=~ |\tilde{g}(\tilde{\nabla}\tilde{r},\nu)|\frac{\tilde{r}}{t}
\sqrt{(1+|Df_{q(t)}|^2)}}\\
&\geq & |\tilde{g}(\tilde{\nabla}\tilde{r},\nu)|
\frac{\tilde{r}}{t}~=~ |\tilde{g}(\tilde{\nabla}\tilde{r},\nu)|
\sqrt{1+ \frac{|\gamma(t)-f(p)|^2}{t^2}}
\geq |\tilde{g}(\tilde{\nabla}\tilde{r},\nu)|.
\end{eqnarray*}
By using  l'H\^{o}pital's Rule,~
$\lim_{t\to +\infty}\frac{\gamma(t)}{t}=\lim_{t\to +\infty}\gamma'(t)$,
and we conclude that
$\lim_{t\to +\infty}|\tilde{g}(\tilde{\nabla}\tilde{r},\nu)|=0$.
\qed\\[5mm]
Set for each $r>0$, $V_M(r)=V_M(\tilde{B}_r(x)\cap M)$. 
\begin{lemma} Let $C_n$ and $ F_n$ be positive constants such that
for $r$ sufficiently large,
$$  F_nr^n \leq V_M(r)\leq C_nr^n.$$
Let $\tau>0$ be a small constant and, $0<a_m<b_m<d_m$ be constants such that
$$ \frac{F_n}{C_n}\left(\frac{b_m}{d_m}\right)^n-
\left(\frac{a_m}{d_m}\right)^n\geq \tau>0.$$
Then,   for any $\epsilon>0$, ~
$\frac{V_M(\frac{b_m}{\epsilon})-V_M(\frac{a_m}{\epsilon})}
{V_M(\frac{d_m}{\epsilon})}\geq \tau. $
\end{lemma}
\noindent
Note that
it is sufficient to take $b_m$ close to $d_m$ and $a_m$
sufficiently smaller than $b_m$ to obtain such a lower bound $\tau$.
 For example, set
$$\theta^n=\frac{C_n}{F_n}\geq 1, \quad
b_m=\frac{d_m}{2}, \quad  a_m=\frac{d_m}{(2^{\alpha}\theta)}<b_m,$$
where $\alpha>1$ and $ \tau=\theta^{-n}(2^{-n}-2^{-\alpha n})>0$.\\[3mm]
\em Proof. \em
$$\frac{V_M(\frac{b_m}{\epsilon})-V_M(\frac{a_m}{\epsilon})}{
V_M(\frac{d_m}{\epsilon})}\geq \frac{F_n}{C_n}\left(\frac{b_m}{d_m}\right)^{n}
-\left(\frac{a_m}{d_m}\right)^{n}\qed$$
Now we follow a similar construction as in J.Li  \cite{Li}.
For each $m$ we consider positive constants 
$0<c_m<a_m<b_m< d_m$ that we
will define later, and a smooth function $\psi_m(t)$ such that $|\psi_m|\leq 1$ and 
$$\psi_m(t)=\left\{\begin{array}{ll}
1 & \mbox{if}~~~a_m\leq t\leq b_m\\
0 & \mbox{if}~~~t\leq c_m~\mbox{or}~t>d_m
\end{array}\right.$$
with
$|\psi'_m(t)|\leq C_m$, and $|\psi''_m(t)|\leq C_m,$
where $C_m>0$ is a positive constant that depends on $a_m$, $b_m$, $c_m$
and $d_m$, that can be given of the form
$$ C_m= E(\frac{1}{a_m-c_m}+\frac{1}{d_m-b_m}
+\frac{1}{(a_m-c_m)^2}+\frac{1}{(d_m-b_m)^2}),$$
where $E>0$ is a constant.
We fix a decreasing sequence $\epsilon_k\to 0$, of positive reals,
and take
$$\eta_{k,m}:= \frac{1}{V_M\left(D\La{(}\frac{c_m}{\epsilon_k},
\frac{d_m}{\epsilon_k}\La{)}\right)},$$
where $D(c,d)=\{y\in M: c\leq \tilde{r}(y)\leq d\}$. 
Now we fix $\lambda>0$ and  consider the function
$$u_{k,m}(y)=\sqrt{\eta_{k,m}}\psi_m(\epsilon_k\tilde{r}(y))e^{i\sqrt{\lambda}
\tilde{r}(y)}.$$
We have
$$\nabla u_{k,m}=
\sqrt{\eta_{k,m}}\La{(}\epsilon_k\psi'_m(\epsilon_k\tilde{r}(y))
+i\sqrt{\lambda}\psi_m(\epsilon_k\tilde{r}(y))\La{)}e^{i\sqrt{\lambda}
\tilde{r}(y)} \nabla \tilde{r}.$$
Set $\zeta^2=1-\xi^2$. Then
\begin{eqnarray}
\lefteqn{\Delta u_{k,m}+\zeta^2\lambda u_{k,m}=}\nonumber\\
&=&\label{l11}\epsilon_{k} \sqrt{\eta_{k,m}}
e^{i\sqrt{\lambda}\tilde{r}}\La{(}\epsilon_k \psi''_m(\epsilon_k\tilde{r})
+2i \sqrt{\lambda}\psi'_m(\epsilon_k\tilde{r})\La{)}
\|\nabla \tilde{r}\|^2\\
&&\label{l12}+ \sqrt{\eta_{k,m}}
e^{i\sqrt{\lambda}\tilde{r}}\La{(}\epsilon_k \psi'_m(\epsilon_k\tilde{r})
+i \sqrt{\lambda}\psi_m(\epsilon_k\tilde{r})\La{)}
\Delta \tilde{r}\\
&&\label{l13} + \lambda (\zeta^2-\|\nabla \tilde{r}\|^2) u_{k,m}.
\end{eqnarray}
Now we prove that~(\ref{l11}) and~(\ref{l12}) tend to $0$ in $L^2$, when
we chose $a_m,b_m,c_m,d_m$ and $\epsilon_k$ in a suitable way 
and let $m,k$ go to $+\infty$.
\begin{eqnarray}
\int_M|(\ref{l11})|^2dV_M&\leq &\eta_{k,m}\epsilon_k^2\int_{D(\frac{c_m}{\epsilon_k},
\frac{d_m}{\epsilon_k})}(2\epsilon_k^2C_m^2+4\lambda{C_m^2})dV_M\nonumber\\
&\leq& \epsilon_k^2{2C^2_m}
\left({\epsilon_k^2}+ 2\lambda\right),
\end{eqnarray}
and using Lemma~\ref{lemma1}
\begin{eqnarray}
\int_M|~(\ref{l12})|^2dV_M&\leq &\eta_{k,m}\int_{D(\frac{c_m}{\epsilon_k},
\frac{d_m}{\epsilon_k})}2(\epsilon_k^2C_m^2+\lambda)
\frac{n}{\tilde{r}^2}dV_M\nonumber\\
&\leq & \frac{\epsilon_k^2}{c_m^2} 2(\epsilon_k^2C_m^2+\lambda)n.
\end{eqnarray}
Now we consider~(\ref{l13}). We have
$\zeta^2-\|\nabla \tilde{r}\|^2_{g_M} 
= -\xi^2+\tilde{g}(\tilde{\nabla}\tilde{r},\nu)^2$.
Then, 
\begin{eqnarray}
\int_M|~(\ref{l13})|^2dV_M&\leq &\int_{D(\frac{c_m}{\epsilon_k},
\frac{d_m}{\epsilon_k})}\lambda^2\eta_{k,m}
(\zeta^2-\|\nabla \tilde{r}\|^2)^2 \nonumber\\
&\leq& \lambda^2 
\sup_{D(\frac{c_m}{\epsilon_k},
\frac{d_m}{\epsilon_k})}(\xi^2-|\tilde{g}(\tilde{\nabla}\tilde{r},\nu)|^2)^2.
\end{eqnarray}
Thus, for any choice of  $a_m,b_m,c_m,d_m\to +\infty$,
 and taking $k=k_m\to +\infty$,
 such that $\epsilon_{k_m}C_m\to 0$,
we may assume $(16)\to 0$, 
and so, we obtain a sequence $u_{k_m,m}\in L^2$ 
s.t. $\Delta u_{k_m,m}+\zeta^2\lambda
u_{k_m,m}\to 0$.
On the other hand, we have
\begin{eqnarray*}
\int_M u_{k,m}^2dV_M &\geq& \eta_{k,m}\int_{D(\frac{a_m}{\epsilon_k},
\frac{b_m}{\epsilon_k})}\!\!\!\!dV_M=\frac{V_M(D(\frac{a_m}{\epsilon_k},
\frac{b_m}{\epsilon_k}))}{V_M(D(\frac{c_m}{\epsilon_k},
\frac{d_m}{\epsilon_k}))}\geq \frac{V_M(\frac{b_m}{\epsilon_k})-
V_M(\frac{a_m}{\epsilon_k})}{V_M(\frac{d_m}{\epsilon_k})}.
\end{eqnarray*}
Choosing $a_m, b_m,c_m,d_m$ satisfying Lemma 3.2,
we have $\int_M u_{k,m}^2\geq \tau>0$, and we may chose
$a_m, b_m,c_m,d_m\to +\infty$ to obtain a sequence $u_{k_m,m}$ that
spans an infinite dimensional subspace of $L^2$. Thus $\zeta^2\lambda$
belongs to the essential spectrum, for any $\lambda>0$.
This proves Theorem~\ref{thm1}.\qed
\begin{remark}{\rm  
Both in Lemma 2.4 and Theorem 2.2 of \cite{Chen}
it is only required 
that  $\lim_{t\to +\infty}\sup_{\tilde{r}(F(y))\geq t}
\tilde{r}(F(y))\|A(y)\|=0$ holds in order
to obtain~(\ref{cond1}) with $C_n=\kappa(M)\omega_n$ and 
$\lim_{t\to +\infty}\inf_{\tilde{r}\geq t}\|\nabla \tilde{r}\|^2=1$. 
We note that
our assumption in the first limit in Proposition 1.1 is equivalent
to the above one, for, if such $\lim$ on $\sup_{\tilde{r}\geq t}
\tilde{r}\|A\|$ 
is zero, then $\lim_{\tilde{r}\to +\infty}\tilde{r}\|A\|$ 
exists and  is zero as well, and vice versa.
Similarly,
  $\lim_{t\to +\infty}\inf_{\tilde{r}\geq t}\|\nabla \tilde{r}\|=1$ 
is equivalent to 
$\lim_{\tilde{r}\to +\infty}\|\nabla \tilde{r}\|^2=1$, 
for $\|\nabla\tilde{r}\|\leq 1$.
The later is equivalent to (\ref{cond2}) with $\xi=0$.}\end{remark}

The next proposition is  obtained from an argument used in the proof of 
Theorem 3.1 of \cite{Do2}, of which  we give here a proof for the sake of 
completeness:
\begin{proposition} If $M$ is a complete noncompact Riemannian manifold
with $\mathrm{Ricci}^M\to -(n-1)c$, when $r(y)\to +\infty$, where $c\geq 0$
is a constant, and $r(y)$ is the intrinsic distance in $M$ to 
a fixed point $x\in M$, then $\lambda_{ess}\leq\frac{(n-1)^2}{4}c$.
\end{proposition}
\noindent
\em Proof. \em By assumption, for any $\delta>0$, $\exists r_0>0$,
such that $\mathrm{Ricci}^M+(n-1)c\geq -\delta$,
for all $y\in M$ with $r(y)\geq r_0$.
Fix a sequence $r_i\to + \infty$. We can find $y_1\in M$ sufficiently far away
from $x$ such that on $B_{r_1}(y_1)$, $\mathrm{Ricci}^M\geq-(n-1)
( c+\frac{1}{2})$.
Next, we  take $y_2$ sufficiently far away from $x$ and such that
$B_{r_2}(y_2)\subset M\backslash B^M_{r_1}(y_1)$ and on $B_{r_2}(y_2)$,
$\mathrm{Ricci}^M\geq -(n-1)(c+\frac{1}{2^2})$. 
By induction we construct a sequence
$y_i$, and balls $B_{r_i}(y_i)\subset M\backslash \cup_{s=1}^{i-1}B_{r_s}
(y_s)$ where $\mathrm{Ricci}^M\geq -(n-1)(c+\frac{1}{2^i})$. From Cheng's
 eigenvalue comparison inequality, on each ball $B_{r_i}(y_i)$, 
$\lambda_1(r_i):=\lambda_1(B_{r_i}(y_i))\leq \lambda_1(D_{r_i})$, where 
$D_{r_i}$ is the disk of radius $r_i$ of the $n$-dimensional space form of
constant sectional curvature $-c_i=-(c+\frac{1}{2^i})$.
Then $\lambda_1(r_i)\leq (n-1)^2 c_i/4+ \psi(r_i)$
where $\psi(r_i)\to 0$ when $r_i\to + \infty$. We consider 
$u_i\in L^2(M)$ the solution of the
Dirichlet problem $\Delta u_i+\lambda_1(r_i)u_i=0$, $u_i=0$ on 
$\partial B_{r_i}(y_i)$
( $u_i$ extended to zero
 on $M\backslash B_{r_i}(y_i)$), 
and $\int_Mu_i^2=1$.
 Let $\lambda_1$ be an accumulation point of
$\lambda_1(r_i)$. Then for a subsequence, $\Delta u_i+\lambda_1u_i\to 0$
in $L^2$, what shows that $\lambda_1 \in \sigma_{ess}(M).$
 This proves  that
$\lambda_{ess}\leq (n-1)^2c/4$.\qed
\begin{remark}{\rm \label{rem4}
In the case of a graphic minimal hypersurface, we have
$$\begin{array}{l}
A(X,Y)= (0,D^2f(X,Y))^{\bot}=\frac{1}{\sqrt{1+|Df|^2}}D^2f(X,Y)\nu,\\[3mm]
Ricci^M(X,X)= -\sum_i|A(X,e_i)|^2\leq 0, \quad \quad
s^M=-\|A\|^2, \end{array}$$
and $\Delta f=0$, that is, $\Delta^0 f-2\frac{|Df|^4}{(1+|Df|^4)^2}
D^2f({Df},{Df})=0$, where $\Delta^0$ is the Euclidean Laplacian. 

We recall the example of
Bombieri, de Giorgi and Giusti \cite{BdGG},
 of non-linear minimal graphic hypersurface
$\Gamma_f$, with $f: \R^{2m}\to\R$, $m\geq 4$, 
is obtained as a limit of solutions
$f^R$ defined 
on balls $B_R$  satisfying the minimal surface equation with
$f^R=f_1$ on $\partial B_R$, and it satisfies
$ |f_1(x)|\leq |f(x)|\leq |f_2(x)|$ where $f_1$ and $f_2$ are defined as follows. Set
$\alpha = (2p+1-\sqrt{\delta})/4>1$,
$\delta= 4p^2-12p +1$,  $p=m-1$, $\lambda$ such that,
$\alpha((2p+1)/(2p+2))<\lambda<\min\{\alpha, p/\alpha^2\}$,
 $D=D(\lambda, p)$, $B=B(\lambda,p)$  sufficiently large
positive constants,
 and
$$P(z)=\int_0^zexp\left( B\int_{|w|}^{\infty}t^{\lambda-2}
(1+ t^{2\alpha(\lambda-1)})^{-1}dt\right)dw.$$
Then,
$$ \begin{array}{l}
 f_1(x)=(u^2-v^2)(u^2+v^2)^{\alpha-1}\\[2mm]
 f_2= P\left((u^2-v^2)+f_1[1+ D|(u^2-v^2)/(u^2+v^2)|^{\lambda-1}]\right),
\end{array}$$
where  $u=(x_1^2+\ldots+x_m)^{1/2}$, and 
$v=(x_{m+1}^2+\ldots+x_{2m}^2)^{1/2}$. 
For $x\in B_h(0)$ the following gradient estimate holds:
$$|D f(x)|\leq c_1 exp\La{(} c_2 \frac{1}{2h}
\sup_{B_{2h}(0)}|f_2|\La{)}=C(h),$$
and one has 
$$\lim\sup_{|x|\to +\infty} \frac{|f(x)|}{|x|^{2\alpha}}
\geq \lim\sup_{|x|\to +\infty} \frac{|f_1(x)|}{|x|^{2\alpha}}=1.$$
This does not allow us to conclude anything about
$\lim_{\tilde{r}\to +\infty}d\tilde{r}(\nu)$. This is an unknown limit
as can been seen from
the behaviour of the second fundamental form
of $\Gamma_f$  given in remark 5  in section 4.
The study of the spectrum of such examples seems to require
the understanding of the behaviour of $d\tilde{r}(\nu^R)$
at $B_R$ of the  solution
$f^R$, using
methods from partial differential equations, a study that is outside the
scope of this paper.}
\end{remark}
\section{Proof of Theorem~\ref{thm2}}
We recall the following lemma of Escobar-Freire (\cite{EF}, Lemma~3.1)
(here we change the sign of $\lambda$ and denote their $f$ by $h$):
\begin{lemma}(\cite{EF}) Let $D$ be a bounded domain with $C^2$ boundary 
in a Riemannian
manifold $M$. Let $u$ and $h$ be functions in $C^1(\bar{D})$. Then
for any $\lambda\in \R$
\begin{eqnarray*}
\int_{D}(\|\nabla u\|^2-\lambda u^2)\Delta h
-2\int_DHess\, h(\nabla u,\nabla u) -2\int_D (\Delta u+ \lambda u)
g_M(\nabla h,\nabla u)\\
= \int_{\partial D}(\|\nabla u\|^2 -\lambda u^2)\frac{\partial h}{
\partial n}-2\int_{\partial D}g_M(\nabla h, \nabla u)\frac{\partial u}{
\partial n}.\quad\quad\quad
\end{eqnarray*}
\end{lemma}
We are considering $F:M\to \mathbb{R}^{n+1}$ a complete minimal
properly immersed hypersurface. Now we follow \cite{EF} closely,
but use extrinsic distance instead of intrinsic.
We take the function $h=\frac{1}{2}\tilde{r}^2$  restricted to $M$.
Then $\|\nabla h\|\leq \tilde{r}$ and $\Delta  h=n$.
We assume $A$ satisfies the conditions of Theorem~\ref{thm2}. Then by (3), 
$Hess\, h(X,X)$ is bounded,   and for any  $X\in T_yM$,
\begin{equation}
 \mathrm{Hess}\, h(X,X)\geq 0, 
\end{equation}
with $ \mathrm{Hess}\, h_{y_0}>0$ at some point $y_0\in M$.
Let us assume that $\lambda$ is an eigenvalue.
We take  $u\in \mathcal{D}(\Delta)$ non zero, such that
$\Delta u +\lambda u=0$.  By the unique continuation
property, $u$ (or $\nabla u$) cannot identically vanish in any open set.
On each extrinsic ball $\tilde{B}_r$ ( that means
$\tilde{B}_r(x)\cap M$), we have
$$\int_{\tilde{B}_r}\|\nabla u\|^2-\lambda u^2
=\int_{\tilde{B}_r}\|\nabla u\|^2+ u\Delta u=
\int_{\tilde{B}_r}\ha \Delta u^2=\int_{\partial \tilde{B}_r}
\frac{1}{2}\frac{\partial u^2}{\partial n}.$$
Then applying Lemma 4.1 and the above equality, we have
\begin{eqnarray}
\lefteqn{\int_{\tilde{B}_r} \mathrm{Hess}\, h(\nabla u, \nabla u)=}
\nonumber\\
&=&\frac{n}{2}\int_{\partial \tilde{B}_r}
\frac{\partial u^2}{\partial n}-\int_{\partial \tilde{B}_r}
(\|\nabla u\|^2-\lambda u^2)\frac{\partial h}{\partial n}
+ 2
g_M(\nabla h, \nabla u)\frac{\partial u}{\partial n}.
\end{eqnarray}
Since $u$ and $\nabla u$ are both in $L^2(M)$, and
$\|\nabla\tilde{r}\|\leq 1$
then, by the co-area formula $$\int_{\epsilon}^{+\infty}dt\int_{
\partial \tilde{B}_t}(\|\nabla u\|^2 +u^2)dS
=\int_{M\backslash \tilde{B}_{\epsilon}}\|\nabla \tilde{r}\|
(\|\nabla u\|^2 +u^2)<+ \infty.$$ 
\begin{lemma}(\cite{EF}) If $\varphi:[\epsilon, +\infty)\to [0,+\infty)$  is a 
mensurable function, with $\epsilon>0$, and such that 
$\int_{\epsilon}^{+\infty}\varphi(t)dt<+\infty$, then there
exists $t_i\nearrow +\infty$, such that $t_i\varphi(t_i)\to 0$.
\end{lemma}
\noindent
\em Proof. \em If such sequence $t_i$ did not exist, then
there exists a constant $C>0$  and a $t_0\geq \epsilon$ such that
for all $t\geq t_0$, $t\varphi(t)\geq C$. But then $\int_{t_0}^{+\infty}
\phi(t)\geq C\log t]^{\infty}_{t_0}=+\infty$, contradicting the
assumption.\qed\\[3mm]
Applying the lemma to $\varphi(t)=\int_{
\partial \tilde{B}_t}(\|\nabla u\|^2 +u^2)dS$, we conclude 
for a sequence $r_i\nearrow +\infty$, ~
$r_i\int_{
\partial \tilde{B}_{r_i}}(\|\nabla u\|^2 +u^2)dS \to 0$.
Thus, (18)  with $r=r_i$ converges to zero when $r_i\to +\infty$ 
( note that $\|\frac{\partial h}{\partial n}\|\leq \|\nabla h\|\leq \tilde{r}$
and $\frac{\partial u^2}{\partial n}=2u \frac{\partial u}{\partial n}$).
Therefore we arrive at
$$\int_{M} \mathrm{Hess}\, h(\nabla u, \nabla u) dV_M=0.$$
Convexity of $h$ (17) implies $ \mathrm{Hess}\,h(\nabla u, \nabla u)=0$ 
everywhere.
Since $h$ is strictly convex at $y_0$, $\nabla u$ must vanish
on a neigbourhood of a point, yielding that $u$
vanishes everywhere. This finishes the proof of
Theorem~\ref{thm2}. The corollaries are an immediate
consequence of Theorem~\ref{thm2} and Remark~\ref{rem4}. \qed
\begin{remark}{\rm   As mentioned in  \cite{Chen},  the examples
of minimal graphs (hence with one end) of Bombieri,
de Giorgi and Giusti for $n\geq 8$ satisfy
$$ +\infty> \sup_r 
\frac{ V_M(M\cap \tilde{B}_r(x))}{r^n} > \omega_n. $$
By Proposition 1.1,
this means that $\lim_{ \tilde{r}(F(y))\to + \infty}
 \tilde{r}(F(y))\|A_y\|=0$ does not hold. Furthermore, we do not know if
(\ref{cond2}) is satisfied and, as a consequence, we cannot exclude
the existence of minimal graphs
with no trivial spectrum.
}\end{remark}
\begin{remark}{\rm  We cannot replace the condition 
$\|A(X,X)\|\leq 1/\tilde{r}$
by the weaker one $\|A\|\leq c/\tilde{r}$, $c>0$ a constant,
 if we want to use the argument 
of Escobar-Freire. Indeed, if we try to 
 rescale $F$ by taking $\hat{F}=cF$,
$\hat{F}$  defines an immersion $\hat{M}$ with
the induced metric
$\hat{g}=c^2g_M$. Identify $T_{\hat{F}(y)}\hat{F}(M)$ with
$T_{F(y)}F(M)$, by identifying $T_{cF(y)}\R^{n+1}$ with $T_{F(y)}\R^{n+1}$,
and 
$\hat{\nu}_{\hat{F}(x)}=\nu_{F(x)}$,
$M$ and $\hat{M}$ have the same Levi-Civita connections and 
$\hat{B}_{\hat{F}(y)}(\tilde{X},\tilde{Y})=\frac{1}{c}
B_{F(y)}(\tilde{X},\tilde{Y}),$
where $\tilde{X}_{cF(y)}=\tilde{X}_{F(y)}$ is a tangent vector
(under the above identifications). Moreover,  
$\hat{\Delta}=\frac{1}{c^2}\Delta^M$,
$\hat{H}=\frac{1}{c^2}H=0$, 
and $\hat{\lambda}(M)=c^{-2}\lambda (M)$, $\hat{\lambda}_{ess}=
c^{-2}\lambda_{ess}$, ${\sigma}(\hat{M})=c^{-2}\sigma(M)$ and similar
for the essential  and the pure point spectrum.  
Unfortunately $\tilde{r}(\hat{F}(y))=c\tilde{r}(F(y))$
which means that we have again
$\|\hat{A}(\tilde{X},\tilde{X})\|\leq c/\tilde{r}(\hat{F}(y))$.}\\
\end{remark}
\section*{Acknowledgements} The first author would like to thank the
Erwin Schr\"{o}dinger Institute in Vienna where part of this work was
carried out, for hospitality and support within the scope of the program
{\it Selected Topics in Spectral Theory}.
The second author would like to thank the hospitality
of professor Francisco Mart\'{\i}n, and 
of the Department of Geometry and Topology of the University of Granada,
where part of this work was completed.

\end{document}